\documentclass[oneside,english]{amsart}
\usepackage[T1]{fontenc}
\usepackage[latin1]{inputenc}
\usepackage{amssymb}
\usepackage{mathrsfs, tabularx}
\usepackage[all]{xypic}
\makeatletter
\theoremstyle{plain}
\newtheorem{thm}{Theorem}[section]
\theoremstyle{plain}
\newtheorem{prop}[thm]{Proposition}
\theoremstyle{plain}

\theoremstyle{remark}
\newtheorem*{rem*}{Remark}
\theoremstyle{definition}
\newtheorem{definition}[thm]{Definition}
\theoremstyle{remark}
\newtheorem{rem}[thm]{Remark}
\newtheorem*{acknowledgement*}{Acknowledgements}

\usepackage{babel}
\makeatother
\begin{document}

\title[A complete lift for semisprays]
{A complete lift for semisprays}
\author[Bucataru]{Ioan Bucataru}
\address{Ioan Bucataru, Faculty of Mathematics, Al.I.Cuza
University, B-dul Carol 11, Iasi, 700506, Romania}
\urladdr{http://www.math.uaic.ro/\textasciitilde{}bucataru/}

\author[Dahl]{Matias F. Dahl}
\address{Matias F. Dahl, Institute of Mathematics, P.O.Box 1100, 02015
Helsinki University of Technology, Finland}

\date{\today}

\begin{abstract} In this paper, we define a complete lift for semisprays. If $S$ is a
semispray on a manifold $M$, its complete lift is a new semispray
$S^c$ on $TM$. The motivation for this lift is two-fold: First,
geodesics for $S^c$ correspond to the Jacobi fields for $S$, and
second, this complete lift generalizes and unifies previously
known complete lifts for Riemannian metrics, affine connections,
and regular Lagrangians. When $S$ is a spray, we prove that the
projective geometry of $S^c$ uniquely determines $S$.  We also
study how symmetries and constants of motions for $S$ lift into
symmetries and constants of motions for $S^c$.

\end{abstract}

\maketitle



\section{Introduction}

Vertical and complete lifts for geometric objects from the base
manifold to its tangent bundle have been studied in various
contexts by many researchers. These lifts have been used to study
various aspects of differential geometry and its applications:
Jacobi fields in \cite{yano73}, stochastic Jacobi fields in
\cite{arnaudon98}, harmonicity in \cite{oniciuc00}, symplectic
geometry in \cite{plebanski01}, symmetries and constants of motion
in \cite{crampin86}. These lifts were introduced by Yano and
Kobayashi in \cite{yano66}, and for a detailed presentation we
refer to the book of Yano and Ishihara \cite{yano73}.

It is known that for a Riemannian metric on a manifold $M$ its
complete lift is a semi-Riemannian metric on $TM$, whose geodesics
are Jacobi fields on $M$, \cite{yano73}. This result has been
extended by Casciaro and Francaviglia in \cite{casciaro97},
N\'u\~nez-Y\'epez and Salas-Brito in \cite{nunez00} and Delgado et
al. in \cite{delgado04}. In these papers, the authors show that if
a system of second order differential equations (SODE) can be
written as a variational problem, then the corresponding system of
Jacobi equations can also be written as a variational problem. For
the homogeneous case, Michor \cite{michor96} showed that one can
lift a spray into a vector field whose integral curves project to
the Jacobi fields of the given spray. This lifted vector field,
which is not a spray, was also studied by Lewis in \cite{lewis98}
and called the tangent lift. For the affine case, Lewis has also
shown how to modify this lift to preserve sprays, \cite{lewis98}.

The main goal of this paper is to generalize and unify the above
lifts into a complete lift for semisprays. Therefore, we will
start with a semispray $S$, whose geodesics are determined by a
system of SODE, define its complete lift, $S^c$, and study its
geodesics. In Theorem \ref{thm_jacobi} we prove that the complete
lift $S^c$ of a semispray $S$ is again a semispray and the
geodesics of $S^c$ are Jacobi fields of $S$. In other words, we
show that the system of geodesic equations of a semispray and the
corresponding system of Jacobi equations are related by the
complete lift. In Section 4.3 we show how this complete lift
contains as special cases the complete lift for a semi-Riemannian
metric \cite{yano66}, an affine connection \cite{lewis98, yano66},
and a regular Lagrangian \cite{casciaro97, delgado04, nunez00}. In
Theorem \ref{thm_jacobi} we also prove that any Jacobi field on a
compact interval can be obtained from a geodesic variation. This
result is well known in the Riemann-Finsler context, \cite{bao00,
docarmo92, sakai92}, where the proof is based on the exponential
map. For the general case of a semispray, the exponential map
cannot be defined, due to the lack of homogeneity. We prove the
correspondence between Jacobi fields and geodesic variations of a
semispray by using the geodesic flow of $S^c$, that is by studying
integral curves of $S^c$ in the second order iterated tangent
bundle $TTM$.

The complete and vertical lifts for geometric objects from a
manifold $M$ to its tangent bundle $TM$ are well known. However,
to lift a semispray we need to lift objects from $TM$ to the
second iterated tangent bundle $TTM$. Since $TTM$ has two vector
bundle structures over $TM$, there are two ways to define such
lifts. The lifting process in this paper has the important feature
that it preserves important geometric objects; semi-Riemannian
metrics, regular Lagrangians, homogeneous functions and vector
fields, semisprays and sprays are all preserved. To avoid studying
separately the vertical and complete lifts from $M$ to $TM$, from
$TM$ to $TTM$, and so on, we introduce in Section 3 a unifying
lifting process from $T^rM$ to $T^{r+1}M$, for any $r\geq 0$. For
$r=0$, this reduces to the vertical and complete lifts for
geometric objects from $M$ to $TM$, as defined by Yano and
Kobayashi \cite{yano66}.

For a spray, the geodesics, up to orientation preserving
reparameterizations, uniquely determine the projective class of a
spray but not the spray itself, \cite{crampin07, douglas27,
matveev05, shen01}. For example, in a subset of the plane, the
Euclidean, Funk and Hilbert metrics all have straight lines as
geodesics, but with possibly different parameterizations,
\cite{shen01}. The main result of Section 5 is Theorem
\ref{thm_projective}. It shows that the projective class of $S^c$
uniquely determines $S$. That is, if sprays $S_1$ and $S_2$ have
the same Jacobi fields, up to orientation preserving
reparameterizations, then $S_1=S_2$.

In the last section of the paper we use the complete lift to study
symmetries and constants of motion of the system of Jacobi
equations in terms of the symmetries and constants of motion of
the geodesic equations of a given semispray. We prove that the
vertical and complete lifts for a constant of motion (or symmetry)
of a system of geodesics equations are constants of motion (or
symmetries) for the corresponding system of Jacobi equations. This
way, we generalize known results for Jacobi fields from Riemannian
geometry, \cite {docarmo92}. In the Lagrangian context, it follows
that the constant of motion considered by N\'u\~nez-Y\'epez and
Salas-Brito in \cite{nunez00} for the system of Jacobi equations
is the complete lift of the Lagrangian function. Also, the
symmetries considered by Arizmendi et al. in \cite{arizmendi03}
for the system of Jacobi equations can be obtained by taking the
complete lift of symmetries of Euler-Lagrange equations. In the
Hamiltonian context, constants of motions for the variational
equations of a given Hamiltonian system have been studied by Case
in \cite{case85} and Giachetta et al. in \cite{giacheta03}.

\section{Higher order iterated tangent bundle}

We start with a real, $C^{\infty}$-smooth, $n$-dimensional
manifold $M$. By $(TM, \pi_0, M)$ we denote its first order
tangent bundle. Every coordinate chart $(U, (x^i))$ on $M$ induces
a coordinate chart $(\pi_0^{-1}(U), (x^i, y^i))$ on $TM$. For
$r\geq 1$, we will denote by $(T^rM, \pi_{r-1}, T^{r-1}M)$ the
$r$th order iterated tangent bundle. Let us also denote $T^0M=M$.
For $r\in \{2,3\}$, the induced local coordinates on $TTM$ and
$TTTM$ will be denoted by $(x,y,X,Y)$ and $(x,y,X,Y,u,v, U, V)$,
respectively.

For a manifold $N$, we denote by $C^{\infty}(N)$ the set of smooth
real functions on $N$, and by ${\mathfrak X}(N)$ the set of smooth
vector fields on $N$. Throughout the paper, all objects will be
assumed to be smooth.

For $r\geq 1$, the slashed tangent bundles $T^rM\setminus\{0\}$
are open subsets in $T^rM$ defined as follows. For $r=1$,
\begin{eqnarray*}
TM\setminus\{0\}=\{(x,y)\in TM, y\neq 0\}, \end{eqnarray*} and for
$r\geq 2$,
\begin{eqnarray*}
\label{slashed_TM} T^rM\setminus\{0\}=\{\xi \in T^rM,
D\pi_{T^{r-1}M\to M}(\xi) \in TM\setminus \{0\} \}.
\end{eqnarray*} It follows that the slashed tangent bundles are
preserved by the tangent map of the canonical submersions. That
is, for $r\geq 2$, $D\pi_{r-2}\left(T^rM\setminus\{0\}\right)
=T^{r-1}M\setminus\{0\}$. For $r=2$ we have explicitly
\begin{eqnarray*}
TTM\setminus\{0\}=\{(x,y,X,Y)\in TTM, X\neq 0\}, \end{eqnarray*}
and $D\pi_{0}\left(TTM\setminus\{0\}\right) =TM\setminus\{0\}$.

For $r\geq 2$, $k\geq 1$, and a map $c: (-\varepsilon,
\varepsilon)^k \to  T^rM$, $\left(t^1,...,t^k\right) \mapsto
c\left(t^1,...,t^k\right)=\left(x^i\left(t^1,...,t^k\right)\right)$,
we define the \emph{derivative} of $c$ with respect to the
variable $t^j$ as the map $\partial_{t^j}c: (-\varepsilon,
\varepsilon)^k \to T^{r+1}M$ defined by $\partial_{t^j}c =
\left(x^i, \partial x^i/\partial t^j\right)$. When $k=1$ we also
write $c'=\partial_tc$.

For $r \geq 0$, a curve $c: I \to T^{r}M$ is said to be a
\emph{regular curve} if $c^{\prime}(t)\in T^{r+1}M\setminus\{0\}$
for all $t$ in $I$, where $I\subset {\mathbb R}$ is an open
interval. For a regular curve $c: I \to T^{r}M$, we have that
$c^{\prime}: I \to T^{r+1}M$ (for $r\geq 0$) and $\pi_{r-1}\circ
c: I \to T^{r-1}M$ (for $r\geq 1$) are also regular curves.

For $r\geq 2$, the $r$th iterated tangent bundle possesses at
least two different vector bundle structures over $T^{r-1}M$. One
is given by the canonical projection $\pi_{r-1}: T^rM
\longrightarrow T^{r-1}M$ and another one is given by the
projection $D\pi_{r-2}: T^rM \longrightarrow T^{r-1}M$. The
involution map $\kappa_r$, interchanges these two vector bundle
structures of $T^rM$.

\begin{definition} For $r\geq 2$, the \textit{involution map} $\kappa_r$, is the unique
diffeomorphism $\kappa_r: T^rM \longrightarrow T^rM$ that
satisfies $\partial_t\partial_s c(t,s)=\kappa_r
\partial_s\partial_t c(t,s)$, for any map $c: (-\varepsilon,
\varepsilon)^2 \longrightarrow T^{r-2}M$. For $r=1$, we define
$\kappa_1=\operatorname{id}_{TM}$.
\end{definition}

It follows immediately that $\kappa_r$ is an involution, which means
that $\kappa_r^2=\operatorname{id}_{T^rM}$. Moreover, $\kappa_r:
(T^rM, \pi_{r-1}, T^{r-1}M) \to (T^rM, D\pi_{r-2}, T^{r-1}M)$ is an
isomorphism of vector bundles. Hence, the following diagram is
commutative.
\begin{eqnarray*}
\begin{xy}
\xymatrix@C=2.0pc @R=3.0pc{
  T^rM\ar@{<->}[rr]^{\kappa_r} &  & T^rM\ar[dl]^{D\pi_{r-2}} \\
& T^{r-1}M \ar@{<-}[ul]^{\pi_{r-1}}& }
\end{xy}
\end{eqnarray*}

In local coordinates, involution maps $\kappa_2: TTM \longrightarrow
TTM$ and $\kappa_3: TTTM \longrightarrow TTTM$ are given by
\begin{eqnarray*}
\kappa_2(x,y,X,Y) & =& (x,X,y,Y), \\
\kappa_3(x,y,X,Y,u,v,U,V)& = &(x,y,u,v,X,Y,U,V).
\end{eqnarray*}

For $r\geq 2$, we have the identities
\begin{eqnarray}
\label{kappar}
  \pi_r\circ D\kappa_r &=& \kappa_r\circ \pi_r,\\
  \label{com_rules}
  D\pi_{r-2} &=&  \pi_{r-1} \circ \kappa_r, \\
\label{kappaId3}
  DD\pi_{r-2} \circ \kappa_{r+1} &=& \kappa_r \circ DD\pi_{r-2}.
\end{eqnarray}
Further properties of involution maps $\kappa_2$, $\kappa_3$, and
so on, are presented in \cite{bd08}. For a discussion about
$\kappa_2$ and the two vector bundle structures of $TTM$ see
\cite{besse78, kolar93, michor00}. The notion of a double vector
bundle structure was introduced by Pradines in \cite{pradines74}.
In \cite{konieczna99}, Konieczna and Urba\'nski have shown that
the framework of double vector bundles is suitable for concepts
such as vertical and complete lifts, linear connections, and
Poisson and symplectic structures.

\section{Vertical and complete lifts}

The vertical and complete lifts for geometric objects from a base
manifold to its tangent bundle has been introduced by Yano and
Kobayashi in \cite{yano66} and studied by numerous geometers in
various settings, see \cite{arnaudon98, oniciuc00, paternain99,
plebanski01, yano73}. The lifts from $T^rM$ to $T^{r+1}M$, studied
in this paper, are motivated by the following aspects. For $r=0$,
we recover the usual vertical and complete lifts, \cite{yano66}.
For $r=1$, the complete lift coincides with the complete lift for
a geodesic spray, introduced by Lewis \cite{lewis98}. For $r\geq
1$, the complete lift preserves important geometric objects:
regular Lagrangians, homogeneous functions and vector fields,
semisprays and sprays.

\subsection{Vertical and complete lifts for functions}

\begin{definition} Let $f\in C^{\infty}(T^rM)$, for some $r\geq 0$.
Then, the \textit{vertical lift} of $f$ is the function $f^v\in
C^{\infty}(T^{r+1}M)$, defined by
\begin{eqnarray} \label{def_fv} f^v\left(\xi\right)=\left(f\circ \pi_r\circ
\kappa_{r+1}\right)\left(\xi\right), \quad \forall \xi \in
T^{r+1}M,
\end{eqnarray} and the \textit{complete lift} of $f$ is the
function $f^c\in C^{\infty}(T^{r+1}M)$, defined by
\begin{eqnarray} \label{def_fc} f^c\left(\xi\right)=
df\left(\kappa_{r+1}\left(\xi\right)\right), \quad \forall \xi \in T^{r+1}M.
\end{eqnarray}
\end{definition}
For $r=0$, formulae \eqref{def_fv} and \eqref{def_fc} coincide
with the usual vertical and complete lift of a function $f\in
C^{\infty}(M)$, \cite{yano66}.

For $r\geq 1$, identity \eqref{com_rules}, implies that the
vertical lift of a function $f\in C^{\infty}(T^rM)$ can be
expressed as follows: $f^v=f\circ D\pi_{r-1}$.

For $r=1$, we have that for $f\in C^{\infty}(TM)$, the vertical and
complete lifts $f^v, f^c\in C^{\infty}(TTM)$ have the following
expressions in local coordinates:
\begin{eqnarray}
\label{vertical_f}
f^v(x,y,X,Y)&=&f(x,X), \vspace{2mm}\\
f^c(x,y,X,Y)&=&\frac{\partial f}{\partial x^i}(x,X)y^i +
\frac{\partial f}{\partial y^i}(x,X)Y^i.\label{complete_f}
\end{eqnarray}

\begin{rem} \label{rem_cfunction} The vertical and complete lifts also
generalize to functions defined on a slashed tangent bundle. If
$f\in C^{\infty}\left( T^rM\setminus\{0\}\right)$, for some $r\geq
1$, then equations \eqref{def_fv} and \eqref{def_fc} define lifts
$f^v, f^c$ and one can prove that $f^v, f^c \in C^{\infty}\left(
T^{r+1}M\setminus\{0\}\right)$.
\end{rem}

\subsection{Complete lift for Lagrangians} \label{sec_Lag}

We show now that the complete lift for functions preserves regular
Lagrangians and contains as a special case the complete lift for
semi-Riemannian metrics introduced by Yano and Kobayashi
\cite{yano66}. Consider $L\in C^{\infty}\left(TM\right)$ a regular
Lagrangian, which means that the Hessian of $L$,
\begin{equation}
g_{ij}(x,y)=\frac{1}{2}\frac{\partial^2 L}{\partial y^i\partial
y^j}(x,y), \label{metricL}
\end{equation} with respect to the fibre coordinates $y$,
has maximal rank $n$ on $TM$, \cite{abraham78, krupkova97}.
According to equation \eqref{complete_f}, the complete lift for a
regular Lagrangian $L$ is given by the following formula
\begin{equation}
L^c(x,y,X,Y)=\frac{\partial L}{\partial x^i}(x,X)y^i +
\frac{\partial L}{\partial y^i}(x,X)Y^i. \label{complete_L}
\end{equation}
It follows that the Hessian of $L^c\in C^{\infty}(TTM)$, with
respect to the fibre coordinates $X, Y$ is given by the following
matrix
\begin{equation}
 \left( \begin{array}{cc} \left(g_{ij}\right)^c  & \left(g_{ij}\right)^v \\ \left(g_{ij}\right)^v & 0
\end{array} \right), \label{glc}
\end{equation}
which has rank $2n$. Therefore, $L^c$ is a regular Lagrangian.
From its local expression \eqref{complete_L} it follows that $L^c$
coincides with the \emph{first-order deformed Lagrangian}
$\mathcal{L}_1$ introduced and studied by Casciaro et al. in
\cite{casciaro96} and Casciaro and Francaviglia in
\cite{casciaro97}. Also the complete lift $L^c$ coincides with the
Lagrangian $\gamma$ considered by N\'u\~nez-Y\'epez and
Salas-Brito in \cite{nunez00} and Delgado et al. in
\cite{delgado04}.

Consider now the particular case, when the regular Lagrangian
$L_g(x,y)=g_x(y,y)$ is induced by a semi-Riemannian metric $g$. By
formula \eqref{complete_L} we have
\begin{equation}
\left(L_g\right)^c(x,y,X,Y)=\frac{\partial g_{ij}}{\partial
x^k}y^kX^iX^j + 2g_{ij}(x)X^iY^j=g^c_{(x,y)}\left((X,Y),
(X,Y)\right). \label{Lcg}
\end{equation}
For the last equality in formula \eqref{Lcg} we used the
components of the complete lift $g^c$ of a semi-Riemannian metric
$g=(g_{ij})$, defined by Yano and Kobayashi \cite{yano66}
\begin{equation}
g^c  =  \left( \begin{array}{cc} \displaystyle\frac{\partial
g_{ij}}{\partial x^k}y^k & g_{ij} \\ g_{ij} & 0 \end{array}
\right). \label{gc}
\end{equation}
Formula \eqref{Lcg} can be written as follows:
\begin{equation}
\left(L_g\right)^c=L_{g^c}, \label{Lcg2}
\end{equation}
and expresses the compatibility between the complete lift
$\left(L_g\right)^c$, of the regular Lagrangian $L_g$ and the
complete lift $g^c$ of the semi-Riemannian metric $g$.

\subsection{Vertical and complete lifts for vector fields}

\begin{definition} \label{definitionAv} For $r\geq 0$, the
\emph{vertical lift} of a vector field $A\in {\mathfrak X}(T^rM)$ is
the vector field $A^v\in {\mathfrak X}(T^{r+1}M)$, defined by
\begin{equation}
A^v(\xi) = D\kappa_{r+1}\circ \partial_s \left. \left(
\kappa_{r+1}(\xi) + s A\circ \pi_r\circ
\kappa_{r+1}(\xi)\right)\right|_{s=0},  \quad \xi\in T^{r+1}M.
\label{defAv}
\end{equation}
\end{definition}
Identity \eqref{kappar} shows that $ \pi_{r+1}\circ A^v =
\operatorname{id}_{T^{r+1}M}$, so $A^v$ is a vector field.

For $r=0$, Definition \ref{definitionAv} reduces to the usual
definition of the vertical lift \cite{abraham78, lewis98, yano66}.

For $r=1$ and a vector field $A\in {\mathfrak X}(TM)$, with local
expression
\begin{eqnarray}
A= \left(x,y,A^i, B^i\right) = A^i(x,y)\frac{\partial}{\partial x^i}
+ B^i(x,y) \frac{\partial}{\partial y^i}, \label{vectorA}
\end{eqnarray}
its vertical lift $A^v\in {\mathfrak X}(TTM)$ is given by
\begin{eqnarray}
A^v =\left(x,y,X,Y, 0, (A^i)^v, 0, (B^i)^v\right) =
\left(A^i\right)^v \frac{\partial}{\partial y^i} +
\left(B^i\right)^v \frac{\partial}{\partial Y^i}. \label{verticalA}
\end{eqnarray}

\begin{definition} \label{definitionAc}
For $r\geq 0$, the \emph{complete lift} of a vector field $A \in
{\mathfrak X}(T^rM)$ is the vector field $A^c \in {\mathfrak
X}(T^{r+1}M)$, defined by
\begin{equation}
A^c = D\kappa_{r+1} \circ \kappa_{r+2} \circ DA \circ
\kappa_{r+1}. \label{defAc}
\end{equation} \end{definition}
Identities \eqref{kappar} and \eqref{com_rules} show that $
\pi_{r+1}\circ A^c = \operatorname{id}_{T^{r+1}M}$, so $A^c$ is a
vector field.

For $r=0$, formula \eqref{defAc} coincides with the usual complete
lift for a vector field $A\in {\mathfrak X}(M)$, \cite{yano66}.

For $r=1$, formula \eqref{defAc} was also consider by Lewis in
\cite{lewis98} to lift an affine spray from $TM$ to $TTM$. For $r=1$
and the vector field $A\in {\mathfrak X}(TM)$, locally given by
formula \eqref{vectorA}, its complete lift $A^c\in {\mathfrak
X}(TTM)$ is given by
\begin{eqnarray}
\nonumber
A^c&=& \left(x,y,X,Y, (A^i)^v, (A^i)^c, (B^i)^v, (B^i)^c\right) \\
   &=& \left(A^i\right)^v \frac{\partial}{\partial x^i} +
\left(A^i\right)^c \frac{\partial}{\partial y^i} +
\left(B^i\right)^v \frac{\partial}{\partial X^i} +
\left(B^i\right)^c \frac{\partial}{\partial Y^i}. \label{completeA}
\end{eqnarray}

\begin{rem} \label{rem_cvector} The vertical and complete lifts also generalize to vector fields defined
on a slashed tangent bundle. If $A\in {\mathfrak X}\left(
T^rM\setminus\{0\}\right)$, for some $r\geq 1$, then equations
\eqref{defAv} and \eqref{defAc} define lifts $A^v, A^c$ and one
can prove that $A^v, A^c \in {\mathfrak X}\left(
T^{r+1}M\setminus\{0\}\right)$.
\end{rem}

The vertical and complete lifts defined in this section have
similar properties as the corresponding lifts studied by Yano and
Ishihara in \cite{yano73}, for $r=0$. For $r\geq 0$,  $f,g\in
C^{\infty}(T^rM)$, and $A,B\in {\mathfrak X}(T^rM)$ we have the
following formulae:
\begin{eqnarray}
\label{fgcv} (fg)^v=f^vg^v, \quad (fg)^c=f^vg^c+ f^cg^v; \\
\label{abcv} [A^c,B^c]=[A,B]^c, \quad [A^c, B^v]=[A,B]^v, \quad
[A^v, B^v]=0; \\
\label{facv} (fA)^c=f^cA^v+f^vA^c, \quad (fA)^v=f^vA^v; \\
\label{afcv} (Af)^c=A^cf^c, \quad (Af)^v=A^cf^v=A^vf^c, \quad A^v
f^v=0.
\end{eqnarray}

The vertical and complete lifts introduced in this sections can be
extended to one-forms and tensors on $T^rM$, by requiring that the
lifts are compatible with the tensor contraction. These lifts
generalize the case $r=0$, studied by Yano and Ishihara in
\cite{yano73}.

\subsection{Lifts for vector fields and their flows}

For $r\geq 0$ and a vector field $A\in {\mathfrak
X}\left(T^rM\right)$, consider its flow $\phi: {\mathscr D}(A) \to
T^rM$, where ${\mathscr D}(A)$ is the maximal domain. Then
${\mathscr D}(A)$ is an open set in $T^rM\times {\mathbb R}$,
\cite{abraham78}. For $\xi \in T^rM$, let $\phi_{\xi}: I(\xi) \to
T^rM$ be the integral curve of $A$ such that $\phi_{\xi}(0)=\xi$
and domain $I(\xi)\subset {\mathbb R}$ is maximal.

\begin{thm} \label{completeLiftFlow} For $r\geq 0$, consider a vector field $A\in
{\mathfrak X}\left(T^rM\right)$ and its complete lift $A^c\in
{\mathfrak X}\left(T^{r+1}M\right)$. Suppose that
$$\phi: {\mathscr D}(A) \to T^rM, \quad \phi^c: {\mathscr D}(A^c) \to
T^{r+1}M,$$ are the flows of $A$ and $A^c$, respectively, with
domains
$$ {\mathscr D}(A)\subset T^rM\times {\mathbb R}, \quad {\mathscr
D}(A^c)\subset T^{r+1}M \times {\mathbb R}.$$ Then,
\begin{eqnarray}
\label{domainc} \left(\pi_r\circ\kappa_{r+1}\times
\operatorname{id}_{\mathbb R} \right){\mathscr D}(A^c)={\mathscr
D}(A)
\end{eqnarray}
and
\begin{eqnarray}
\label{flowc} \phi^c_t(\xi) =
  \kappa_{r+1}\circ D\phi_t \circ \kappa_{r+1} (\xi),\quad \left(\xi,t\right) \in {\mathscr
  D}(A^c),
\end{eqnarray}
where $D\phi_t$ is the tangent map of the map $\xi\mapsto
\phi_t(\xi)$, for a fixed $t$.
\end{thm}

\begin{proof}
To prove inclusion ``$\subset$'' in \eqref{domainc} we start with
an integral curve of $A^c$ and show that it projects into an
integral curve of $A$. If $\left(\xi_0, t_0\right)\in {\mathscr
D}(A^c)$, then $\phi^c_{\xi_0}: I(\xi_0) \to T^{r+1}M$ is an
integral curve of vector $A^c$. Using formula \eqref{defAc},
identity \eqref{com_rules} and the commutation rule
$\pi_{r+1}\circ DA= A \circ \pi_r$ we obtain
$$ D\left(\pi_r\circ \kappa_{r+1}\right) \circ A^c= A\circ
\left(\pi_r\circ \kappa_{r+1}\right). $$ It follows that
$\pi_r\circ\kappa_{r+1}\circ \phi^c_{\xi_0}$ is an integral curve
of $A$ defined on $I(\xi_0)$. Therefore
$\left(\pi_r\circ\kappa_{r+1}(\xi_0), t_0\right)\in {\mathscr
D}(A)$ and inclusion ``$\subset$'' in \eqref{domainc} follows.

To prove inclusion ``$\supset$'' in \eqref{domainc} and formula
\eqref{flowc} we start with an integral curve of $A$ and show that
it can be lifted into an integral curve of $A^c$. Let
$\left(\pi_r\circ\kappa_{r+1}(\xi_0), t_0\right)\in {\mathscr
D}(A)$, for some $\xi_0 \in T^{r+1}M$. Then
$\phi_{\pi_r\circ\kappa_{r+1}(\xi_0)}:
I(\pi_r\circ\kappa_{r+1}(\xi_0)) \to T^rM$ is an integral curve of
vector field $A$. For each $t\in
I(\pi_r\circ\kappa_{r+1}(\xi_0))$, consider
\begin{eqnarray}
  \gamma(t) = \kappa_{r+1}\circ D\phi_t \circ \kappa_{r+1}
  (\xi_0).\label{gamma}
\end{eqnarray}
We will prove now that $\gamma$ is well defined and smooth on
$I(\pi_r\circ\kappa_{r+1}(\xi_0))$, $\gamma(0)=\xi_0$, and that
$\gamma$ is an integral curve of $A^c$.

Since ${\mathscr D}(A)$ is open in $T^rM\times {\mathbb R}$, for
any $(\pi_r\circ\kappa_{r+1}(\xi_0),t)\in {\mathscr D}(A)$ there
exists a neighborhood $U\subset T^rM$ of
$\pi_r\circ\kappa_{r+1}(\xi_0)$ and a neighborhood $I\subset
I(\pi_r\circ\kappa_{r+1}(\xi)) \subset {\mathbb R}$  of $t$ such
that $U\times I\subset {\mathscr D}(A)$. For vector
$\kappa_{r+1}(\xi_0)\in T_{\pi_r\circ\kappa_{r+1}(\xi_0)}U$, we
can find a smooth curve $v\colon (-\varepsilon, \varepsilon) \to
U$ such that $v(0)=\pi_r\circ\kappa_{r+1}(\xi_0)$ and
$\kappa_{r+1}(\xi_0)=\partial_sv(0)$. It follows that
$\phi(v(s),\tau)$ is well defined and smooth on the open
neighborhood $(-\varepsilon, \varepsilon)\times I$ of $\{0\}\times
\{t\}$. Therefore, for all $t$ in
$I(\pi_r\circ\kappa_{r+1}(\xi_0))$
\begin{eqnarray} \label{var} \gamma(t) = \kappa_{r+1}\circ \partial_s\phi(v(s),\tau)|_{s=0, \tau=t},\end{eqnarray} and hence curve
$\gamma$ is well defined and smooth on
$I(\pi_r\circ\kappa_{r+1}(\xi_0))$. Setting $t=0$ in equation
\eqref{gamma}, and using $\phi_0 = \operatorname{id}_{T^{r}M}$ we
have that $\gamma(0)=\xi_0$. By formula \eqref{var} we have that
\begin{eqnarray*} \label{acgamma}
  (A^c\circ \gamma)(t) &=& D\kappa_{r+1} \circ \kappa_{r+2} \circ DA \circ D\phi_t \circ
  \kappa_{r+1}(\xi)\\
             &=& D\kappa_{r+1} \circ \kappa_{r+2} \circ \partial_s \left(A \circ \phi(v(s), \tau)\right) |_{s=0, \tau=t}\\
             &=& D\kappa_{r+1} \circ \kappa_{r+2} \circ \partial_s \partial_{\tau} \phi(v(s),\tau) |_{s=0, \tau=t} \\
             &=& D\kappa_{r+1} \circ \partial_{\tau} \partial_s \phi(v(s),\tau) |_{s=0, \tau=t} \\
             &=& D\kappa_{r+1} \circ \partial_{\tau} \left(\phi_{\tau} \circ\kappa_{r+1}(\xi_0)\right)|_{\tau=t} \\
             &=& D\kappa_{r+1} \circ \partial_{\tau} \left(\kappa_{r+1} \circ \gamma\right)(\tau)|_{\tau=t} \\
             &=& \gamma'(t).
\end{eqnarray*}
In the above calculations we used that $A\circ \phi_{\tau} =
\partial_{\tau} \phi_{\tau}$ and the definitions of $\kappa_{r+2}$ and
$\gamma$.  The last equality follows since $ D\kappa_{r+1}\circ
\partial_{\tau}\gamma =  \partial_{\tau}\kappa_{r+1}(\gamma).$

We have shown that $\gamma$ is an integral curve of $A^c$ defined
on $I(\pi_r\circ\kappa_{r+1}(\xi_0))$, such that
$\gamma(0)=\xi_0$. Therefore, $(\xi_0,t_0)\in {\mathscr D}(A^c)$,
which implies that $(\pi_r\circ\kappa_{r+1}(\xi_0),t_0)\in
\left(\pi_r\circ\kappa_{r+1}\times \operatorname{id}_{\mathbb R}
\right){\mathscr D}(A^c)$ and inclusion ``$\supset$'' in
\eqref{domainc} is true. It also follows that formula
\eqref{flowc} is true.
\end{proof}

From formula \eqref{domainc} we have that ${\mathscr
D}(A^c)=T^{r+1}M\times {\mathbb R}$ if and only if ${\mathscr
D}(A)=T^{r}M\times {\mathbb R}$. That is, vector field $A\in
{\mathfrak X}\left(T^rM\right)$ is complete if and only if $A^c\in
{\mathfrak X}\left(T^{r+1}M\right)$ is complete.

\section{Semisprays and their complete lift}

The geometry of a system of second order differential equations
have been initiated in the first half of the last century through
the work of Douglas \cite{douglas27}, Cartan \cite{cartan33},
Chern \cite{chern39}, Kosambi \cite{kosambi33}, and others. Such a
system of SODE on $M$ can be represented using a semispray, which
is a special vector field on the phase space $TM$. For a
semispray, its geodesics are defined as projections onto $M$ of
its integral curves in $TM$. In this section we prove that the
complete lift for a semispray $S$ on $M$ is a semispray $S^c$ on
$TM$, and show how geodesics of $S^c$ are related to Jacobi fields
of $S$.

\subsection{Semisprays and their geodesics}

\begin{definition} \label{def:spray} For $r\geq 1$, a \emph{semispray $S$ on $T^{r-1}M$}
is a vector field $S\in {\mathfrak
X}\left(T^rM\setminus\{0\}\right)$ such that $\kappa_{r+1}\circ
S=S$.
\end{definition}

From the definition, we see that a semispray on $T^{r-1}M$ is a
section for both bundle structures over $T^{r}M$. This can
equivalently be written as follows: a vector field $S\in
{\mathfrak X}\left(T^rM\setminus\{0\}\right)$ is a semispray on
$T^{r-1}M$ if and only if $D\pi_{r-1}\circ
S=\operatorname{id}_{T^rM\setminus\{0\}}$.

For $r=1$, Definition \ref{def:spray} reduces to the usual
definition of a semispray, \cite{bucataru07}. In local coordinates
$(x,y)$, a semispray $S$ on $M$ has the form
\begin{eqnarray}
\nonumber S & = &  \left( x^i,y^i,y^i,-2G^i(x,y) \right) \\
\label{SprayDef1} & = & y^i \frac{\partial}{\partial x^i} -2
G^i(x,y) \frac{\partial}{\partial y^i},
\end{eqnarray}
for some functions $G^i$ defined on the domain of the considered
induced chart on $TM\setminus\{0\}$.

In local coordinates $(x,y, X, Y)$, a semispray $S$ on $TM$ has the
form
\begin{eqnarray}
\nonumber
S &=&  \left( x^i,y^i,X^i, Y^i, X^i, Y^i, -2G^i(x,y,X,Y), -2H^i(x,y,X, Y) \right) \\
\label{SprayDef2} &=& X^i \frac{\partial}{\partial x^i} + Y^i
\frac{\partial}{\partial y^i} -2 G^i(x,y,X,Y)
\frac{\partial}{\partial X^i} - 2 H^i(x,y,X,Y)
\frac{\partial}{\partial Y^i},
\end{eqnarray}
for some functions $G^i, H^i$ defined on the domain of the
considered induced chart on $TTM\setminus\{0\}$.

\begin{definition} \label{def_geod} For $r\geq 1$, a \emph{geodesic} for a
semispray $S$ on $T^{r-1}M$, is a regular curve $c\colon I \to
T^{r-1}M$ such that $c''=S(c')$, and $I\subset \mathbb{R}$ is the
maximal domain.\end{definition}

From formula \eqref{SprayDef1} we see that a regular curve
$c\colon I \to M$, $c(t)=(x^i(t))$, is a geodesic for a semispray
$S$ on $M$ if and only if it locally satisfies the following
system of second order ordinary differential equations:
\begin{eqnarray}
\label{geo}
  \frac{d^2x^i}{dt^2} + 2 G^i\left(x, \frac{dx}{dt}\right) &=& 0.
\end{eqnarray}

We refer to equations \eqref{geo} as to the \textit{geodesic
equations} of the semispray $S$, given by formula \eqref{SprayDef1}.

\subsection{Jacobi fields and geodesic variations}

Let $c: I \to M$ be a geodesic for a semispray $S$ on $M$. A
geodesic variation of $c$ is a smooth map $V: I\times
(-\varepsilon, \varepsilon) \to M$, $V=V(t,s)$ such that
\begin{itemize}
\item[1)] $V(t,0)=c(t)$, for all $t$ in $I$,
\item[2)] $V(t,s)$ is a geodesic for all $s$ in $(-\varepsilon,
\varepsilon)$. \end{itemize} The variation vector field $J(t)=\left.
\partial_s V(t,s)\right|_{s=0}=(x^i(t), y^i(t))$ satisfies the following system of second order
ordinary differential equations:
\begin{eqnarray}
\label{Jacobi}
  \frac{d^2y^i}{dt^2} + 2 \frac{\partial G^i}{\partial y^j}\left(x, \frac{dx}{dt}\right)
  \frac{d y^j}{dt} + 2 \frac{\partial G^i}{\partial x^j}\left(x, \frac{dx}{dt}\right)y^j = 0.
\end{eqnarray}
We refer to equations \eqref{Jacobi} as to the \textit{Jacobi
equations} or the \textit{variational equations} of the system of
geodesic equations \eqref{geo}. Geodesic variations and the
corresponding Jacobi equations \eqref{Jacobi} were considered
by Kosambi in \cite{kosambi33}, Cartan in \cite{cartan33} and Chern
in \cite{chern39} for arbitrary systems of second order differential
equations.

\begin{definition} \label{def_jacobi}
A vector field $J:I \to TM$, $J(t)=\left(x^i(t), y^i(t)\right)$, along a geodesic
$c=\pi_0\circ J$ of a semispray $S$ on $M$, is called a \emph{Jacobi
field} if it satisfies the Jacobi equations \eqref{Jacobi}.
\end{definition}

Using the dynamical covariant derivative induced by $S$, one can
express the geodesic equations \eqref{geo} and the Jacobi
equations \eqref{Jacobi} in an invariant way, \cite{bucataru07,
jerie02}. We will see also, from Theorem \ref{thm_jacobi}, that
the system of geodesic equations and the system of Jacobi
equations for a semispray, together, form the geodesic equations
of the complete lift of the semispray. Hence, the systems of
equations \eqref{geo} and \eqref{Jacobi} do not depend on a
particular choice of local coordinates.

For a geodesic $c$ of a semispray $S$ on $M$, its tangent vector
$c'(t)$ is a Jacobi field along $c$ that is locally determined by
geodesic variation $V(t,s)=c(t+s)$.

We define the \textit{complete lift}, $S^c$, of a semispray $S$ as
the complete lift of $S$ as a vector field, and the
\textit{geodesic flow} of semispray $S$ as the flow of $S$ as a
vector field. The domain ${\mathscr D}(S)$ of the geodesic flow of
$S$ is open in $TM\setminus\{0\}\times {\mathbb R}$ and the domain
${\mathscr D}(S^c)$ of the geodesic flow of $S^c$ is open in
$TTM\setminus\{0\}\times {\mathbb R}$.

\begin{thm} \label{thm_jacobi} Let $S$ be a semispray on $M$.
\begin{enumerate}
\item[i)]
The complete lift $S^c$ is a semispray on $TM$.
\item[ii)] A curve $J\colon I \to TM$ is a geodesic for the semispray $S^c$ if and only if
$J$ is a Jacobi field along $c=\pi_0\circ J$.
\item[iii)] The domains of the geodesic flows of the two semisprays $S$ and $S^c$ are
related by the formula \begin{equation} (D\pi_0\times
\operatorname{id}_{\mathbb R}){\mathscr D}(S^c) = {\mathscr D}(S).
\label{domain_ssemis}
\end{equation}
\item[iv)] For any restriction of a Jacobi field $J$ to a compact interval
$K$, there is a geodesic variation $V: K\times (-\varepsilon,
\varepsilon) \to M$ of $c=\pi_0\circ J$ such that
$J(t)=\partial_sV(t,s)|_{s=0}$, for all $t\in K$.
\end{enumerate}
\end{thm}
\begin{proof}
i) According to Remark \ref{rem_cvector}, formula \eqref{defAc}
can be used also to lift vector fields defined on slashed tangent
bundles. Therefore, for $r=1$, we have that the complete lift of
$S$ is defined as follows:
\begin{eqnarray} \label{sc} S^c=D\kappa_2 \circ \kappa_3 \circ DS
\circ \kappa_2.\end{eqnarray} If $S$ is given in local coordinates
by formula \eqref{SprayDef1}, then using formula
\eqref{completeA}, its complete lift is given by
\begin{eqnarray}
\nonumber
  S^c &=& \left(x^i,y^i,X^i,Y^i,X^i,Y^i,-2(G^i)^v, -2(G^i)^c \right) \\
      &=& X^i \frac{\partial}{\partial x^i} + Y^i \frac{\partial}{\partial y^i}  -
      2 \left(G^i\right)^v \frac{\partial}{\partial X^i} -2 \left(G^i\right)^c
  \frac{\partial}{\partial Y^i}.\label{compSc}
\end{eqnarray}
It follows that $\kappa_3\circ S^c=S^c$ and hence $S^c$ is a
semispray on $TM$.

ii) A curve $J\colon I \to TM$, $J(t)=(x^i(t), y^i(t))$ is a
geodesic for $S^c$ if and only if
\begin{eqnarray}
\label{eq_giv} \frac{d^2x^i}{dt^2} & = &
-2\left(G^i\right)^v\left(x, y, \frac{dx}{dt},
\frac{dy}{dt}\right); \vspace{2mm}\\
\label{eq_gic} \frac{d^2y^i}{dt^2} & = &
-2\left(G^i\right)^c\left(x, y, \frac{dx}{dt},
\frac{dy}{dt}\right).
\end{eqnarray}
By formulae \eqref{vertical_f} and \eqref{complete_f} equations
\eqref{eq_giv} and \eqref{eq_gic} are equivalent to  geodesic
equations \eqref{geo} and Jacobi equations \eqref{Jacobi},
respectively and ii) follows.

iii) Using similar arguments as we used in the proof of Theorem
\ref{completeLiftFlow}, it follows that the solution $J$ of the
Jacobi equations \eqref{Jacobi} with initial condition $J'(0)=\xi
\in TTM\setminus\{0\}$ and the solution $c$ of the geodesic
equations \eqref{geo} with the initial condition
$c'(0)=D\pi_0(\xi)\in TM\setminus\{0\}$ are defined over the same
open interval $I(\xi)=I(D\pi_0(\xi))\subset {\mathbb R}$.
Therefore we obtain that the domains of the flows for the two
semisprays $S$ and $S^c$ are related by the formula
\eqref{domain_ssemis}.

iv) Consider $J: K\to TM$ a Jacobi field and let $c: K \to M$,
$c=\pi_0\circ J$ be the underlying geodesic, where $K\subset I$ is
a compact interval and $I$ is the maximal domain for $J$. Denote
$\xi=J'(0)\in TTM\setminus\{0\}$.

Let $\phi$ be the geodesic flow of semispray $S$. If
$\phi_{D\pi_0(\xi)}: K\to TM$ is the integral curve of $S$ with
$\phi_{D\pi_0(\xi)}(0)=D\pi_0(\xi)$ then $\pi_0\circ
\phi_{D\pi_0(\xi)}=c$. Since $J'$ is an integral curve of $S^c$
then, using a similar argument as we used in the proof of Theorem
\ref{completeLiftFlow}, we have that
$$ J'(t)=\kappa_2\circ D\phi_t \circ \kappa_2(\xi), \quad \forall t\in I.$$

Since ${\mathscr D}(S)$ is open in $TM\setminus\{0\}\times
{\mathbb R}$ and $\xi \in TTM\setminus\{0\}$ is fixed, for
$\left(D\pi_0(\xi), t\right)\in {\mathscr D}(S)$ there exists
$U_t\subset TM\setminus\{0\}$ an open neighborhood of
$D\pi_0(\xi)$ and $I_t\subset {\mathbb R}$ an open interval that
contains $t$ such that $U_t \times I_t \subset {\mathscr D}(S)$.
For the vector $\kappa_2(\xi)\in T_{D\pi_0(\xi)}U_t$, there exists
a differentiable curve $w: (-\varepsilon_t, \varepsilon_t) \to
U_t$ such that $w(0)=D\pi_0(\xi)$ and
$\kappa_2(\xi)=\partial_sw(0)$.

Since $K$ is compact, from the open covering $\{I_t, t\in K\}$,
one can choose a finite covering $\{I_{t_j}, j\in \{1,...,N\}\}$
of it. Define $\varepsilon:=\min\{\varepsilon_{t_j}, j\in
\{1,...,N\}\}$ and $U:=\cap_{j=1}^NU_{t_j}$. It follows that
$\phi\left(w(s),t\right)$ is well defined for all $(t,s)\in K
\times (-\varepsilon, \varepsilon)$. Consider the geodesic
variation $ V: K\times (-\varepsilon, \varepsilon) \to
TM\setminus\{0\}$, given by
\begin{equation} V(t,s)=\pi_0\circ
\phi\left( w(s),t\right). \label{geod_var} \end{equation} We have
that $V(t,0)=\pi_0\circ \phi\left(D\pi_0(\xi),t\right) =c(t)$.
Moreover, since $J(t)=\pi_1\circ J'(t)$ is a geodesic of $S^c$, we
have
\begin{eqnarray*}
J(t) & = & \pi_1 \circ \kappa_2\circ D\phi_t \circ
\partial_sw(s)|_{s=0} \\
& = & D\pi_0 \circ D\phi_t \circ \partial_sw(s)|_{s=0} \\
& = & D\pi_0 \circ \partial_s\phi(w(s),t)|_{s=0} \\
& = & \partial_s\left(\pi_0 \circ \phi( w(s),t)\right)|_{s=0} \\
& = & \partial_sV(t,s)|_{s=0}.
\end{eqnarray*}
Hence, $J(t)$ is the variation vector field of the geodesic
variation $V(t,s)$ given by formula \eqref{geod_var}.
\end{proof}

\subsection{Discussion of results}

We discuss now the results of Theorem \ref{thm_jacobi} in various
contexts: Riemannian, Finslerian, affine, and Lagrangian.

Last part of the Theorem \ref{thm_jacobi} is usually proved in the
Riemannian context, \cite{docarmo92, sakai92}, and the proof is
based on the exponential map. In the Finslerian context, the
result is discussed in the book \cite{bao00}. Since for a
semispray we do not require any homogeneity condition, their
geodesics do not inherit a homogeneity condition that will allow
us to define the exponential map. Therefore the standard proof of
Theorem \ref{thm_jacobi} iv) in Riemann-Finsler context does not
generalize directly to semisprays. Instead, the above proof relies
on the geodesic flow of $S^c$, that is by studying integral curves
of $S^c$ in the second order iterated tangent bundle $TTM$.

Consider $g$ a semi-Riemannian metric on $M$ and denote by $S_g$
the geodesic spray induced by $g$, which is a semispray on $M$.
Consider $g^c$, the complete lift of $g$, given by formula
\eqref{gc}. Let $S_{g^c}$ be the geodesic spray of the
semi-Riemannian metric $g^c$, which is a semispray on $TM$. Then,
it follows that
\begin{equation} \left(S_g\right)^c=S_{g^c}, \label{sgc}
\end{equation} and therefore the complete lift for a semispray
discussed in Theorem \ref{thm_jacobi} contains as a special case,
the complete lift for a semi-Riemannian metric introduced by Yano
and Kobayashi in \cite{yano66}.

Consider an affine spray $S_{\nabla}$, whose function coefficients
are $G^i(x,y)=\gamma^i_{jk}(x)y^jy^k$, where $\gamma^i_{jk}(x)$
are the coefficients of an affine connection $\nabla$ on $M$.
Consider $\nabla^c$, the complete lift of the affine connection
$\nabla$, defined by $\nabla^c_{X^c}Y^c=\left(\nabla_XY\right)^c$,
\cite{yano66}. Let $S_{\nabla^c}$ be the affine spray induced by
the affine connection $\nabla^c$. Then, it can be shown that
\begin{equation} \left(S_{\nabla}\right)^c=S_{\nabla^c}, \label{snablac} \end{equation}
and therefore the complete lift for a semispray discussed in
Theorem \ref{thm_jacobi} contains as a special case, the complete
lift for an affine connection introduced by Yano and Kobayashi in
\cite{yano66}. Moreover, $\left(S_{\nabla}\right)^c$ coincides
with $D\kappa_2\circ S_{\nabla}^T\circ \kappa_2$ considered by
Lewis in \cite{lewis98}. Lewis also has shown that the geodesics
of the complete lift for an affine spray are Jacobi fields for the
given spray.

Consider a regular Lagrangian $L$ on $TM$, and let $S_L$ be the
corresponding Euler-Lagrange vector field, which is a semispray on
$M$. Consider the complete lift $L^c$, given by equation
\eqref{complete_L}, which is a regular Lagrangian on $TTM$ and
$S_{L^c}$ the corresponding Euler-Lagrange vector field, which is
a semispray on $TM$. Then, it can be shown that \begin{equation}
(S_L)^c=S_{L^c}, \label{slc} \end{equation} and therefore the
complete lift for a semispray discussed in Theorem
\ref{thm_jacobi} contains as a special case, the complete lift for
a Lagrange function discussed in Section 3.2. Formula \eqref{slc}
shows that the Euler-Lagrange equations of $L^c$ consist of the
Euler-Lagrange equations of the Lagrangian $L$ and their Jacobi
equations. This result agrees with the results obtained by
Casciaro and Francaviglia in \cite{casciaro97}, Delgado et al. in
\cite{delgado04}, and N\'u\~nez-Y\'epez and Salas-Brito in
\cite{nunez00}.

\section{Sprays and their complete lifts}

In this section we prove that the complete lift studied in the
previous sections preserves Liouville vector fields, homogeneous
functions and vector fields and sprays. Two sprays are
projectively related if they have the same geodesics up to
orientation preserving reparameterizations. The main result of
this section is Theorem \ref{thm_projective}. It shows that
geometric properties for affine, semi-Riemannian or Finsler sprays
can be determined from studying the corresponding projective
properties of the complete lift. The projective geometry of sprays
was initiated by Douglas in \cite{douglas27}. For a modern
presentation of projective connection we refer to
\cite{crampin07}, for projectively related sprays see Szilasi
\cite{szilasi03}, and for the Finslerian case see Shen
\cite{shen01}.

\subsection{Sprays}

The notion of spray, in the affine context, was introduced and
studied by Ambrose et al. in \cite{ambrose60}. A general spray,
which does not reduce to an affine spray, has to be defined on the
slashed tangent bundle, $T^rM\setminus\{0\}$.

The homogeneity of various objects on tangent bundles can be
characterized, using Euler's theorem, in terms of the Liouville
vector field.

For $r\geq 1$, the \emph{Liouville vector field} ${\mathbb C}_r\in
{\mathfrak X}(T^rM)$ is defined as follows:
$$
{\mathbb C}_r(\xi) = \left. \partial_s\left(
\xi+s\xi\right)\right|_{s=0}, \quad \xi\in T^rM.
$$
From its definition we can see that the Liouville vector field
${\mathbb C}_r$ is a globally defined vector field. In local
coordinates we obtain the following formulae for the Liouville
vector fields ${\mathbb C}_1$ and ${\mathbb C}_2$.
\begin{eqnarray}
 {\mathbb C}_1 = (x,y,0,y) = y^i \frac{\partial}{\partial y^i}, \\
 {\mathbb C}_2 = (x,y,X,Y, 0,0, X, Y) = X^i \frac{\partial}{\partial
 X^i} + Y^i \frac{\partial}{\partial
 Y^i}. \label{EulerVectorField}
\end{eqnarray}
Using the definition of the complete lift we find that the Liouville
vector fields are preserved by the complete lifts, which means that
\begin{eqnarray}
\label{EulerIteration}
  {\mathbb C}_r^c = {\mathbb C}_{r+1}, \quad r\geq 1.
\end{eqnarray}
\begin{definition} \label{def_hom} For $r\geq 1$ and a non-negative integer $s$, we say that a function
$f\colon T^rM\setminus\{0\}\to \mathbb{R}$ is \emph{positively
$s$-homogeneous} if for all $\xi\in T^rM\setminus\{0\}$,
$$
f(\lambda \xi) = \lambda^s f(\xi), \quad \lambda>0.
$$ \end{definition}
Throughout the paper, homogeneity refers to positive homogeneity
only. For $r\geq 1$, Euler's theorem for homogeneous functions
implies that a function $f\in C^{\infty}(T^rM\setminus\{0\})$ is
$s$-homogeneous if and only if ${\mathbb C}_r(f)=sf$.

The homogeneity can be extended to other objects that live on
$T^rM\setminus\{0\}$, see \cite[Section 1.5]{bucataru07} for
$r=1$. For $r\geq 1$, a vector field $A\in {\mathfrak
X}(T^rM\setminus\{0\})$ is $s$-homogeneous if $[{\mathbb C}_r,
A]=(s-1)A$. For example, the Liouville vector field ${\mathbb
C}_r$ is $1$-homogeneous.

Next proposition will show that the complete and vertical lifts
preserve the homogeneity for functions and vector fields.

\begin{prop}  Let $r\geq 1$.
\label{prop_hom}
\begin{itemize} \item[i)] If function $f\in C^{\infty}(T^rM\setminus\{0\})$
is $s$-homogeneous, then its vertical and complete lifts $f^v, f^c
\in C^{\infty}(T^{r+1}M\setminus\{0\})$ are $s$-homogeneous
functions.
\item[ii)] If vector field $A\in {\mathfrak X}(T^rM\setminus\{0\})$ is $s$-homogeneous,
then its vertical and complete lifts $A^v, A^c \in {\mathfrak
X}(T^{r+1}M\setminus\{0\})$ are $s$-homogeneous vector fields.
\end{itemize}
\end{prop}
\begin{proof}
If $f: T^rM\setminus\{0\} \longrightarrow \mathbb{R}$ is
$s$-homogeneous, then equations \eqref{EulerIteration} and
\eqref{afcv} imply that ${\mathbb C}_{r+1}(f^c)={\mathbb
C}_r^c(f^c)=\left({\mathbb C}_r(f)\right)^c=sf^c$ and $f^c$ is
$s$-homogeneous. The other claims follow similarly.
\end{proof}

\begin{definition} \label{spray} For $r\geq 1$, a semispray
$S$ on $T^{r-1}M$ is a
\emph{spray} on $T^{r-1}M$ if it is $2$-homogeneous, which means
that $[{\mathbb C}_r,S]=S$.
\end{definition}

For $r=1$, a spray $S$ on $M$ is given by formula
\eqref{SprayDef1}, where the spray coefficients $G^i(x,y)$ are
$2$-homogeneous functions on $TM\setminus\{0\}$. In this case,
Definition \ref{spray} is equivalent to the usual definition of a
spray, studied in various contexts, \cite{ambrose60, bao00,
bucataru07, shen01, szilasi03}.

If for a spray $S$ on $M$, its domain is the whole $TM$, instead
of $TM\setminus\{0\}$, then the homogeneity condition implies that
the functions $G^i(x,y)$ are quadratic in $y$. It follows that
there exist functions $\gamma^i_{jk}(x)$ on the base manifold $M$
such that $2G^i(x,y)=\gamma^i_{jk}(x)y^jy^k$. Functions
$\gamma^i_{jk}(x)$ are local coefficients of an affine connection
$\nabla$ on the base manifold and the spray is said to be an
\textit{affine spray}, \cite{ambrose60}.

Now, we show that the complete lift preserves sprays.

\begin{prop} \label{cor_spray} If $S$ is a spray on $M$,
then its complete lift $S^c$ is a spray on $TM$.
\end{prop}
\begin{proof} It follows by Theorem \ref{thm_jacobi} i) and Proposition \ref{prop_hom}
ii).
\end{proof}

Since vertical and complete lifts preserve homogeneity, it follows
that for a spray $S$ on $M$, with spray coefficients $G^i$, the
spray coefficients $\left(G^i\right)^v, \left(G^i\right)^c$, of
the complete lift $S^c$ are $2$-homogeneous functions with respect
to the fiber coordinates $X^i, Y^i$. One can use the homogeneity
conditions to show that $tc'(t)$ is also a Jacobi field, and
locally it is determined by the geodesic variation
$V_2(t,s)=c(t+ts)$.

Let us point out that the complete lift in Definition
\ref{definitionAc} is closely related to other lifts. In
\cite{lewis98}, Lewis defined the \emph{tangent lift} of an affine
spray $S$ on $M$ as the vector field $S^T=\kappa_3\circ DS$. This
definition also coincides with the \emph{canonical lift} of a
vector field considered by Fisher and Laquer in \cite{fisher99}.
In \cite{michor96}, Michor proved that if $S$ is a spray, then
integral curves of $S^T$ project onto Jacobi fields of $S$. Here,
$S^T$ is a vector field, but it is not a spray. However, Lewis in
\cite{lewis98}, noticed that replacing $S^T$ with $D\kappa_2\circ
S^T \circ \kappa_2$ one obtains a spray, provided that $S$ is an
(affine) spray. This modified lift coincides with $S^c$ in
Definition \ref{definitionAc}. Integral curves of vector field
$S^T$ and spray $S^c$ are closely related. A curve $\gamma\colon I
\to TTM$ is an integral curve of $S^T=\kappa_3\circ DS$ if and
only if $\kappa_2\circ \gamma$ is an integral curve of $S^c$.

\subsection{Projectively related sprays}

\begin{definition} \label{projs} For $r\geq 1$, two sprays $S_1$ and $S_2$ on $T^{r-1}M$ are said to be
\emph{projectively related} if their geodesics coincide up to
orientation preserving reparameterizations.
\end{definition}

For $r=1$, Definition \ref{projs} reduces to the usual definition
of projectively related sprays on $M$, \cite{shen01, szilasi03}.
Projectively related sprays have the same geodesics as point sets.
For symmetric sprays, the converse is also true, see
\cite{crampin07}.

For $r\geq 1$, two sprays $S_1$ and $S_2$ on $T^{r-1}M$ are
projectively related if and only if there exists a $1$-homogeneous
function $P: T^rM\setminus\{0\} \longrightarrow \mathbb{R}$ such
that
\begin{equation}
S_1=S_2+2P{\mathbb C}_r. \label{gpr}
\end{equation}
This characterization can be expressed in terms of the spray
coefficients as follows: two sprays $S_1$ and $S_2$ on $T^{r-1}M$,
with spray coefficients $G^i_1$ and $G^i_2$, are projectively
related if and only if $G^i_2(x,y)=G^i_1(x,y)+P(x,y)y^i$. For
$r=1$ this characterization of projectively related sprays  is
discussed in \cite{shen01}.

Next theorem states that projective geometry of the complete lift
of a spray uniquely determines the spray.

\begin{thm} \label{thm_projective} If $S_1$ and $S_2$ are sprays on $M$ then the following
statements are equivalent:
\begin{itemize}
\item[i)] $S_1$ and $S_2$ coincide, which means that $S_1$ and $S_2$ have the same
geodesics as parameterized curves; \item[ii)] $S_1^c$ and $S_2^c$
are projectively related, which means that $S_1$ and $S_2$ have
the same Jacobi fields up to orientation preserving
reparameterizations.\end{itemize}
\end{thm}
\begin{proof}
Consider that $S^c_1$ and $S^c_2$ are projectively related. These
two sprays are vector fields on $TTM\setminus\{0\}$ and their
spray coefficients, according to formula \eqref{compSc}, are
$\left(G^i_1\right)^v$, $\left(G^i_1\right)^c$ and
$\left(G^i_2\right)^v$, $\left(G^i_2\right)^c$, respectively.
According to formula \eqref{gpr}, it follows that there is a
$1$-homogeneous function $Q: TTM\setminus\{0\} \longrightarrow
\mathbb{R}$ such that
\begin{eqnarray}
\label{givq} \left(G^i_2\right)^v(x,y,X,Y) & = &
\left(G^i_1\right)^v(x,y,X,Y) + Q(x,y,X,Y) X^i, \\
\label{gicq} \left(G^i_2\right)^c(x,y,X,Y) & = &
\left(G^i_1\right)^c (x,y,X,Y) + Q(x,y,X,Y) Y^i.
\end{eqnarray}
Now, equation \eqref{givq} implies that $Q(x,y,X,Y)$ depends only
on $(x,X)$. Hence $Q$ does not depend on variables $y$ and $Y$ and
therefore $P=Q\circ \mathbb{C}_1\in
C^{\infty}\left(TM\setminus\{0\}\right)$ is a $1$-homogeneous
function such that $Q=P^v$. Now equation \eqref{givq} implies
$G^i_2(x,y)  = G^i_1(x,y) + P(x,y) y^i$, and equation \eqref{gicq}
implies that $(Py^i)^c=P^vY^i$. Equation \eqref{fgcv} and
$(y^i)^c=Y^i$ implies that $P^c=0$. Therefore $P$ is locally
constant on $TM\setminus \{0\}$, and by $1$-homogeneity $P=0$.
\end{proof}

\subsection{Discussion of results}

We discuss now the results of Theorem \ref{thm_projective} in two
contexts: affine and Finslerian.

In the book \cite{yano73}, Yano and Ishihara pointed out that the
complete lift does not preserve the projective class of affine
connections. If we apply Theorem \ref{thm_projective} to the case
of affine sprays, we obtain the following result. Two affine
connections on a manifold $M$ coincide if and only if they have
the same Jacobi fields up to orientation preserving
reparameterizations.

Let $F_1$ and $F_2$ be two Finsler functions, which means that
their Lagrangian functions $L_1=F_1^2$ and $L_2=F_2^2$ are
$2$-homogeneous regular Lagrangians that are smooth only on
$TM\setminus \{0\}$, \cite{bucataru07, shen01}. Let $S_1$ and
$S_2$ be the corresponding geodesic sprays of the two Finsler
functions. According to Theorem \ref{thm_projective}, we have that
$S_1=S_2$ if and only if the two Finsler functions have the same
Jacobi fields up to orientation preserving reparameterizations,
which is equivalent to the fact that $L_1^c$ and $L_2^c$ are
projectively related.

\section{Lie symmetries and constants of motion for Jacobi equations}

We have seen that the system of Jacobi equations \eqref{Jacobi}
can be derived from the system of geodesic equations \eqref{geo}
using the complete lift. In this section we study how Lie
symmetries and constant of motions for a system of SODE can be
lifted to Lie symmetries and constant of motions for the
corresponding system of Jacobi equations.

\begin{definition} Consider a
semispray $S$ on $T^{r-1}M$ for some $r\geq 1$. \begin{itemize}
\item[i)] A function $f\in C^{\infty}(T^rM)$ is a \textit{constant
of motion} for $S$ if $S(f)=0$.
\item[ii)] A vector field $A\in {\mathfrak X}(T^{r-1}M)$ is a \textit{Lie
symmetry} for $S$ if $[S, A^c]=0$.
\end{itemize}
\end{definition} For $r=1$, we refer to \cite{crampin86, krupkova97} for
constants of motions and Lie symmetries for geodesic sprays and
Euler-Lagrange vector fields and to \cite{bucataru07} for constants
of motions and Lie symmetries for semisprays.

\begin{prop} \label{prop_symmetry} Consider a semispray $S$ on $M$.
\begin{itemize} \item[i)] If $f\in C^{\infty}(TM)$ is a constant of motion for the
semispray $S$, then its vertical and complete lifts $f^v, f^c\in
C^{\infty}(TTM)$ are constants of motion for the semispray $S^c$.
\item[ii)] If $A\in {\mathfrak X}(M)$ is a Lie symmetry for the
semispray $S$, then its vertical and complete lifts $A^v, A^c\in
{\mathfrak X}(TM)$ are Lie symmetries for the semispray $S^c$.
\item[iii)] Let $\psi$ be the flow of a Lie symmetry $A\in {\mathfrak
X}(M)$. Then its complete lift $\psi^c$ preserves the Jacobi
fields of $S$.
\end{itemize}
\end{prop}
\begin{proof} Using properties \eqref{afcv}, for $r=1$, we have that
$S^c(f^c)=\left(S(f)\right)^c$ and $S^c(f^v)=\left(S(f)\right)^v$.
Hence, if $S(f)=0$ then $S^c(f^c)=0$ and $S^c(f^v)=0$. Similarly,
using properties \eqref{abcv}, for $r=1$, we have that $[S,A]=0$
implies $[S^c, A^c]=0$ and $[S^c, A^v]=0$.

For the last part, consider a Lie symmetry $A\in {\mathfrak X}(M)$
of a semispray $S$ on $M$, which means that $[S,A^c]=0$. Let
$\phi_t$ the geodesic flow and $\psi_s$ the flow of the vector
field $A$. Then $A$ is a Lie symmetry if and only if $\phi_t\circ
\psi_s= \psi_s\circ \phi_t$, which means that the flow of $A$
preserves the solution curves of the system of geodesic equations
\eqref{geo}, \cite{krupkova97}. Using either Proposition
\ref{prop_symmetry} ii) or Theorem \ref{completeLiftFlow} it
follows that the complete lift $\psi^c$, which is the flow of
$A^c$, preserves the Jacobi fields of the given system of SODE.
\end{proof}

Although the proof of Proposition \ref{prop_symmetry} is very
simple, it generalizes results in various contexts: Lagrangian,
Finslerian or semi-Riemannian.

Proposition \ref{prop_symmetry} provides constants of motion and
Lie symmetries for the system of Jacobi equations \eqref{Jacobi}
when we know such constants of motions and Lie symmetries for the
system of geodesic equations \eqref{geo}. This aspect has been
studied by Case in \cite{case85} and by Giachetta et al. in
\cite{giacheta03} for Hamiltonian systems.

Consider $S$ the Euler-Lagrange vector field of a regular
Lagrangian $L:TM\to {\mathbb R}$, which is a semispray on $M$.
Then, the energy function $E_L={\mathbb C}_1(L)-L$ is a constant
of motion for the semispray $S$. From Proposition
\ref{prop_symmetry} i) it follows that
\begin{eqnarray} \label{constantel}
\left(E_L\right)^c=E_{L^c} \quad \textrm{and} \quad
\left(E_L\right)^v
\end{eqnarray} are constants of motion for the system of Jacobi equations
\eqref{Jacobi}. In equation \eqref{constantel}, $E_{L^c}={\mathbb
C}_2(L^c)-L^c$ is the energy of the lifted Lagrangian $L^c$
introduced in Section 3.2. This aspect has been studied by
Arizmendi et al. in \cite{arizmendi03} and by N\'u\~nez-Y\'epez
and Salas-Brito in \cite{nunez00}.

For the case of a Finsler space, the homogeneity of
the Finsler function $F$ implies that its energy is $E_{F^2}=F^2$.
Therefore,
\begin{eqnarray} \label{constantf2} \frac{\partial F^2}{\partial
x^i}\left(x, \frac{dx}{dt}\right) y^i + \frac{\partial
F^2}{\partial y^i}\left(x, \frac{dx}{dt}\right) \frac{dy^i}{dt}
\quad \textrm{and} \quad F^2\left(x,\frac{dx}{dt}\right)
\end{eqnarray} are constants along the Jacobi field
$J(t)=\left(x^i(t), y^i(t)\right)$ of the Finsler space.

The fact that expressions \eqref{constantel} and
\eqref{constantf2} are constants along the Jacobi fields
generalizes known results from Riemannian geometry, see
\cite{docarmo92, sakai92}. For a semi-Riemannian metric $g$
consider the Lagrangian $L(x,y)=g_{ij}(x)y^iy^j$. Accordingly,
from expressions \eqref{constantf2} we obtain that
\begin{eqnarray}
\label{constant} g\left(\nabla y, \frac{dx}{dt}\right) \quad
\textrm{and} \quad L\left(x,\frac{dx}{dt}\right)
\end{eqnarray} are constants along the Jacobi field $J(t)=(x^i(t),
y^i(t))$ of the semi-Riemannian metric $g$. In expression
\eqref{constant}, $\nabla y$ denotes the covariant derivative of
$y$ with respect to the Levi-Civita connection of the
semi-Riemannian metric $g$.

\begin{acknowledgement*}
I.B has been supported by grant ID 398 from the Romanian Ministry of
Education. M.D. has been supported by Academy of Finland Center of
Excellence Programme 213476, the Institute of Mathematics at the
Helsinki University of Technology, and Tekes project MASIT03 --
Inverse Problems and Reliability of Models.
\end{acknowledgement*}

\end{document}